\documentclass[12pt]{article}
\usepackage{latexsym, amsfonts, amssymb, amscd, amsthm}
\usepackage{amsmath, mathrsfs, yhmath}
\usepackage{epsfig, graphicx, subfig}
\usepackage{hyperref}
\newtheorem{theorem}{Theorem}[section]

\newtheorem{lemma}[theorem]{Lemma}

\newtheorem{example}[theorem]{Example}
\newtheorem{remark}[theorem]{Remark}

\numberwithin{equation}{section}
\allowdisplaybreaks[4]

\title{\bf Star-shaped Curves under Gage's Area-preserving Flow and the CSF}

\author{ \bf  Laiyuan Gao~~~~Shicheng Zhang~~~~Yuntao Zhang}

\date{}

\begin{document}
\maketitle
\noindent{\bf Abstract}
Mayer asks a question what closed, embedded and nonconvex initial curves guarantee that Gage's area-preserving flow (GAPF) exists globally.
A folklore conjecture since 2012 says that GAPF evolves smooth, embedded and star-shaped initial curves globally. In this paper, we prove this conjecture by
using Dittberner's singularity analysis theory.
A star-shaped ``flying wing" curve is constructed to show that GAPF may not always preserve the star-shapedness of evolving curves.
This example is also a negative answer to Mantegazza's open problem whether the curve shortening flow (CSF) always preserves the star shape of the evolving curves.
\\\\
\noindent {\bf Keywords} {curve shortening flow, area-preserving flow, star-shaped, flying wing}

\noindent {\bf Mathematics Subject Classification (2000)} {35C44, 35K05, 53A04}\\\\

\section{Introduction}
Since 1980s, the curvature flow of curves arose some interest among experts of geometric analysis. Given a convex curve in the Euclidean plane, Gage and Hamilton
\cite{Gage-1983, Gage-1984, Gage-Hamilton-1986} proved that this curve, under the {\bf curve shortening flow (CSF)}, shrinks to a point with asymptotically circular shape.
Later Grayson \cite{Grayson-1987, Grayson-1989} proved that a smooth, closed and embedded curve becomes a convex one under the CSF. However, there exist dumbbell-shaped evolving surfaces,
which may split into pieces under the mean curvature flow before the singularities occur \cite{Angenent-1992, Mayer-Simonett-2000}. So, Grayson
type theorem is a typical phenomenon in the planar geometry. And this result has several new proofs
\cite{Andrews-Bryan-2011, Gao-Pan-2016, Hamilton-1995, Huisken-1998, Magni-Mantegazza-2014}.

In the same time, experts also consider curvature flows with a non-local term of closed and embedded curves. For example, Gage \cite{Gage-1986} has studied next area-preserving flow
\begin{equation}\label{eq:1.1.202401}
\left\{\begin{array}{l}
\frac{\partial X}{\partial t}(\varphi, t)=(\kappa(\varphi, t)-\frac{2\pi}{L(t)})N(\varphi, t) \ \ \ \text{in} \ \ S^1\times (0, \omega),\\
X(\varphi, 0)= X_0(\varphi) \ \  \ \ \ \  \text{on} \ \ S^1,
\end{array} \right.
\end{equation}
where $X: S^1\times [0, \omega)\rightarrow \mathbb{R}^2$ is a family of smooth, closed and embedded curves with curvature $\kappa=\kappa(\varphi, t)$,
length $L=L(t)$ and the unit inner normal $N(\varphi, t)$. If the initial curve $X_0$ is convex, then Gage \cite{Gage-1986} proved that the evolving curve $X(\cdot, t)$
is smooth on time interval $[0, +\infty)$ and it is deformed into a circle as time tends to infinity.
Gage noticed that simple closed curves may develop self-intersections and singularities under his area-preserving flow.
Later, Mayer \cite{Mayer-2001} and Mayer-Simonett \cite{Mayer-Simonett-2000} proved Gage's observation by numerical experiments.

An open question is what additional conditions of a closed, embedded and nonconvex initial curve $X_0$ guarantee that {\bf Gage's area-preserving flow (GAPF)} exists globally, where
we say a flow is {\bf global}, if the evolving curve $X(\cdot, t)$ is smooth for every $t\in [0, +\infty)$.
This question was proposed in the last sentence of Mayer's paper \cite{Mayer-2001} in 2001.
In 2012, a folklore conjecture said that the star-shapedness condition of the initial curve is a sufficient one, where we call a closed and embedded curve
{\bf star-shaped} if there is a point $P$, called as a {\bf star center}, such that the {\bf support function} with respect to this point, i.e.
$p(\varphi) = -\langle X_0(\varphi)-P, N_{\text{in}}(\varphi)\rangle$, is positive everywhere.
In this paper, we first reconsider GAPF \cite{Gage-1986} and confirm this conjecture. 
\begin{theorem}\label{thm:1.1.202401}
Let $X_0: S^1 \rightarrow \mathbb{R}^2$ be a smooth, embedded and star-shaped curve in the plane. Gage's area-preserving flow (\ref{eq:1.1.202401}), with this
initial curve, exists globally on time interval $[0, +\infty)$ and deforms the evolving curve into a circle as $t \rightarrow +\infty$.
\end{theorem}

Theorem \ref{thm:1.1.202401} is an alternative answer to Mayer's above question. This result is also a partial generalization of the main theorem in our previous paper \cite{Gao-Pan-2023},
where the authors show that a smooth, embedded, centrosymmetric
and star-shaped curve, evolving as an initial curve of GAPF (\ref{eq:1.1.202401}), yields a family of smooth star-shaped curves on time interval
$[0, +\infty)$ and the evolving curve converges to a circle as $t\rightarrow +\infty$.

\begin{remark}\label{rem:1.2.202401}
Like Grayson Theorem, the fact in Theorem \ref{thm:1.1.202401} is a typical phenomenon in the planar geometry. There exist star and dumbbell shaped evolving surfaces
which may split into pieces under Huisken's volume-preserving mean curvature flow \cite{Angenent-1992, Huisken-1987, Mayer-Simonett-2000} as singularities occur.
In addition, the embedded condition of $X_0$ in Theorem \ref{thm:1.1.202401} can not be dropped either.
See an example of immersed star-shaped curves in Remark 3.4 of the paper \cite{Gao-Pan-2023}, where the immersed and star-shaped evolving
curve under GAPF (\ref{eq:1.1.202401}) blows up in finite time.
\end{remark}

In the year 2010, Mantegazza \cite{Mantegazza-2010} asked whether the CSF preserves the star shape of the evolving curve. One may
propose same question for GAPF. To our surprise, we find a star-shaped curve like a flying wing which may lose its star shape under GAPF (see Figure \ref{fig:4}
and Example \ref{exa:3.1.202401} in Section 3). Motivated by the results in the paper \cite{Gao-Pan-2023}, this kind of star-shaped curves are
far from central symmetry.
In order to answer Mantegazza's question, we prove next result.

\begin{theorem}\label{thm:1.3.202401}
Suppose $X_0$ is a smooth, embedded and star-shaped curve in the plane. Let this curve evolve as an initial curve according to Gage's area-preserving flow (\ref{eq:1.1.202401}).
Then we have a family of smooth and embedded curves in the plane $X: S^1\times [0, + \infty)\rightarrow \mathbb{R}^2$.
Let $A_0$ be the area of the region bounded by $X_0$. If there exists $t_0\in (0, \frac{A_0}{2\pi})$ such that
the curve $X(\cdot, t_0)$ is not star-shaped with respect to any points, then the curve shortening flow with initial $X_0$ also loses the star shape of the evolving curve before
the time $t_0$.
\end{theorem}

Along with above questions proposed by Mayer and Mantegazza, there are in fact many interesting results about curvature flow of star-shaped curves or hypersurfaces,
such as Tsai \cite{Tsai-1996}, Yagisita \cite{Yagisita-2005}, Smoczyk \cite{Smoczyk-1998} and references therein.

This paper is organized as follows. In Section 2, Theorem \ref{thm:1.1.202401} is proved.
In Section 3, we construct a star-shaped curve and explain that GAPF (\ref{eq:1.1.202401}) with this initial curve
may not preserve the star shape of the evolving curve. Theorem \ref{thm:1.3.202401} is proved in Section 4.
All the curves mentioned in this paper are closed and smooth in the plane.

\section{Proof of Theorem \ref{thm:1.1.202401}}
Let $X_0: S^1 \rightarrow \mathbb{R}^2$ be an embedded and smooth curve in the plane. In the year 1996, Chow, Liou and Tsai \cite{Chow-Liou-Tsai-1996} defined a kind of {\bf turning angle} of $X_0$ by
\begin{eqnarray}\label{eq:2.1.202401}
\Theta(X_0) = \inf_{\Gamma} \int_{\Gamma} \kappa ds,
\end{eqnarray}
where the infimum is taken over all connected arcs $\Gamma$ on the curve. They studied an expansion flow of embedded curves with turning angle greater than $-\pi$.
Recently, Dittberner \cite{Dittberner-2021} defines that the curve $X_0$ has {\bf locally total curvature} greater than $-\pi$ if $\int_{\Gamma} \kappa ds > -\pi$
holds for all connected arcs $\Gamma$ on the curve. And he proves the next result.
\begin{lemma}\label{lem:2.1.202401}
Let $X_0$ be a closed and embedded curve with locally total curvature no less than $-\pi$. Gage's area-preserving flow (\ref{eq:1.1.202401}) with this
initial curve exists globally.
\end{lemma}

It follows from previous definitions, an embedded curve has locally total curvature greater than or equal to $-\pi$ if and only if Chow, Liou and Tsai's
turning angle $\Theta(X_0) \geq -\pi$. To avoid ambiguity in other places, we would like to call $\Theta(X_0)$ in (\ref{eq:2.1.202401}) as {\bf Chow-Liou-Tsai turning angle}
of the curve $X_0$.

Lemma \ref{eq:2.1.202401} is a profound result of GAPF. To prove this theorem, Dittberner \cite{Dittberner-2021} first reviews two types of singularities of the flow (\ref{eq:1.1.202401}).
Let $\omega>0$ be the first singular time of the flow (\ref{eq:1.1.202401}). A singularity point $P$ of this flow is called as
Type I, if $\lim\limits_{t \rightarrow \omega-} |\kappa|_{\max}(t)\sqrt{\omega -t}$ is bounded; or called as
Type II, if $\lim\limits_{t \rightarrow \omega-} |\kappa|_{\max}(t) \sqrt{\omega -t}$ is unbounded, where
$$|\kappa|_{\max}(t) = \max\{|\kappa(\varphi, t)|~| \varphi \in S^1\}$$ is the maximum value of the curvature function $\kappa(\varphi, t)$ at time $t$.
To study the singularity model of the flow, one may rescale the curve near the singular time. After a proper rescaling, Dittberner proves that
Type I singularity models are Abresch-Langer curves \cite{Abresch-Langer-1986, Epstein-Weinstein-1987}
and Type II singularity models are the Grim Reapers \cite{Angenent-1991}.
Since the flow (\ref{eq:1.1.202401}) preserves the bounded area of the evolving curve, Type I singularity does not occur.
If the Chow-Liou-Tsai turning angle of $X_0$ is no less than $-\pi$, then the evolving curve preserves this property
along the flow (\ref{eq:1.1.202401}). Dittberner further shows that Huisken's isoperimetric ratio \cite{Huisken-1998}
$$\min_{p, q \in S^1} \frac{d(p, q, t)}{l(p, q, t)}$$
has a positive lower bound which is independent of time. So he could also exclude Type II singularities in any finite time.
Therefore, the flow (\ref{eq:1.1.202401}), under the condition of Lemma \ref{lem:2.1.202401}, exists on time interval $[0, +\infty)$.

\begin{figure}[tbh]
\centering
\includegraphics[scale=0.4]{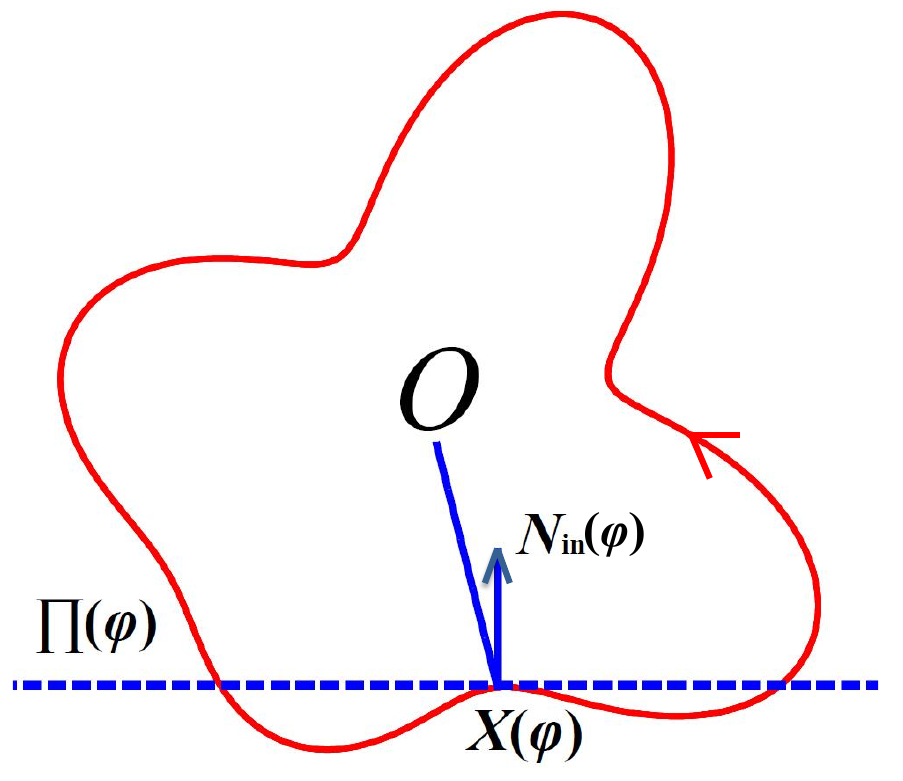}
\caption{A Star-shaped Curve and the Half Plane $\Pi(\varphi)$.}\label{fig:1}
\end{figure}
Our idea of proving Theorem \ref{thm:1.1.202401} is to check that a star-shaped curve has Chow-Liou-Tsai turning angle $\Theta(X_0) > -\pi$.
Let $X_0$ be an embedded and star-shaped $C^1$ curve with a star center $O$. The tangent line at the point $X_0(\varphi)$ divides the plane into two parts.
Denote by $\Pi(\varphi)$ the interior of the half plane that $N_{\text{in}}(\varphi)$ points to. Let $\mathfrak{C}$ be the set of all star centers of the curve $X_0$.
We first prove a geometric property of star-shaped curves.

\begin{lemma}\label{lem:2.2.202401}
The set of all star centers of a star-shaped $C^1$ curve is the intersection of all open half planes pointed by the unit inner normal, i.e.,
\begin{eqnarray}\label{eq:2.2.202401}
\mathfrak{C} = \bigcap_{\varphi \in S^1} \Pi(\varphi).
\end{eqnarray}
\end{lemma}
\begin{proof}
Let the point $P$ be one of the star centers of a star-shaped curve $X_0$. Denote by $p(\cdot)$ the support function with respect to $P$.
Since $p(\varphi)  = - \langle X_0(\varphi)-P, N_{\text{in}}(\varphi)\rangle$ is positive, the star center $P$ lies in the half plane $\Pi(\varphi)$.
By the arbitrary choice of $\varphi$, the point $P$ is in the set $\cap_{\varphi\in S^1} \Pi(\varphi)$.
Since $P$ can be any point in the set of star centers, we have $\mathfrak{C} \subseteq \cap_{\varphi\in S^1} \Pi(\varphi)$.

Once a point $P$ is in the half plane $\Pi(\varphi)$ for a fixed $\varphi \in S^1$, we have the inner product
\begin{eqnarray}\label{eq:2.3.202401}
- \langle X_0(\varphi)-P, N_{\text{in}}(\varphi)\rangle >0.
\end{eqnarray}
For a star-shaped curve $X_0$, we have proved that $\cap_{\varphi\in S^1} \Pi(\varphi)$ is not an empty set in the last paragraph.
Let the point $P$ be in the set $\cap_{\varphi\in S^1} \Pi(\varphi)$.
Then, by (\ref{eq:2.3.202401}), the support function with respect to $P$ is positive everywhere. Therefore, $P$ is a star center.
By the arbitrary choice of $P$, we have the other inclusion $\cap_{\varphi\in S^1} \Pi(\varphi) \subseteq \mathfrak{C}$.
\end{proof}
\begin{figure}[tbh]
\centering
\includegraphics[scale=0.45]{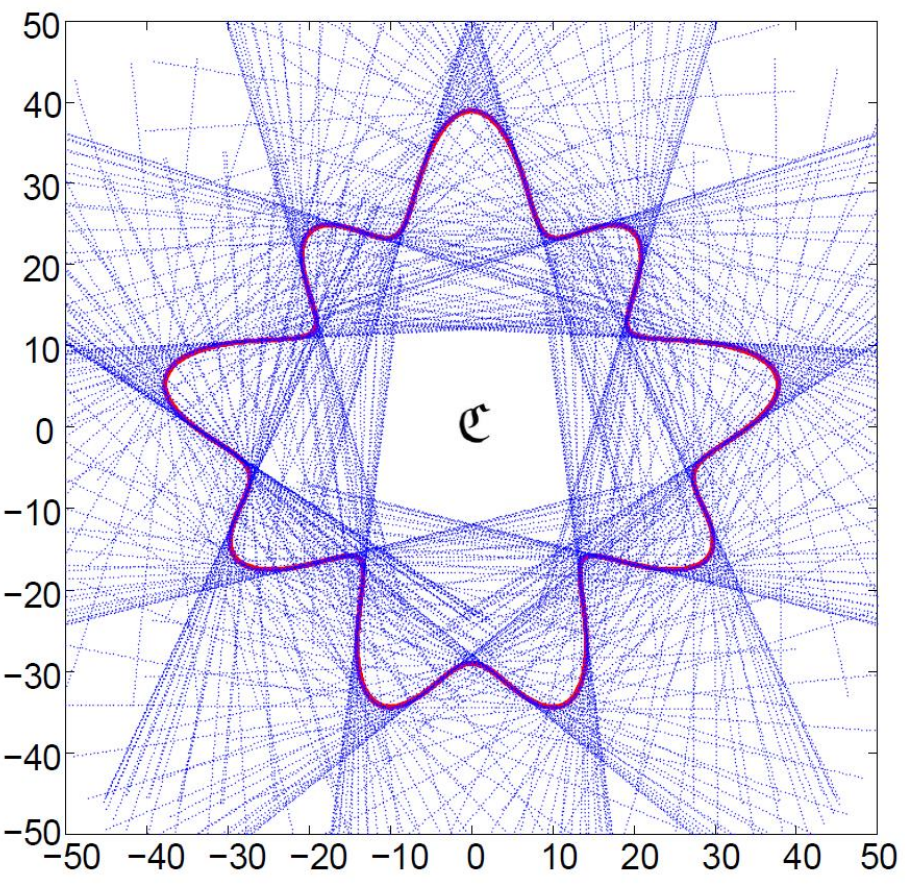}
\caption{A Star-shaped Curve and Its Kernel $\mathfrak{C}$.}\label{fig:2}
\end{figure}
Given a star-shaped curve, we call the set of its all star centers as its {\bf kernel}.
As a corollary of Lemma \ref{lem:2.2.202401}, the kernel of a star-shaped curve is an open and convex set. Smith \cite{Smith-1968} proved this fact by another method in the year 1968 for general star-shaped
sets. An example of this fact is given as follows.

\begin{example}\label{exa:2.3.202401}
Given a smooth and periodic function $r(\theta) = 30+4\cos(4\theta)+5\sin(9\theta)$ on the interval $[0, 2\pi]$. Let $X_0= r(\theta)(\cos \theta, \sin\theta)$ be a curve
with radial function $r(\theta)$. Since $r(\theta)$ is smooth and positive everywhere, the curve $X_0$ is smooth and star-shaped with respect to $O$.
In Figure 2, this curve together its tangent bundle are plotted. Its kernel is bounded by a convex polygon.
\end{example}

\begin{figure}[tbh]
\centering
\includegraphics[scale=0.5]{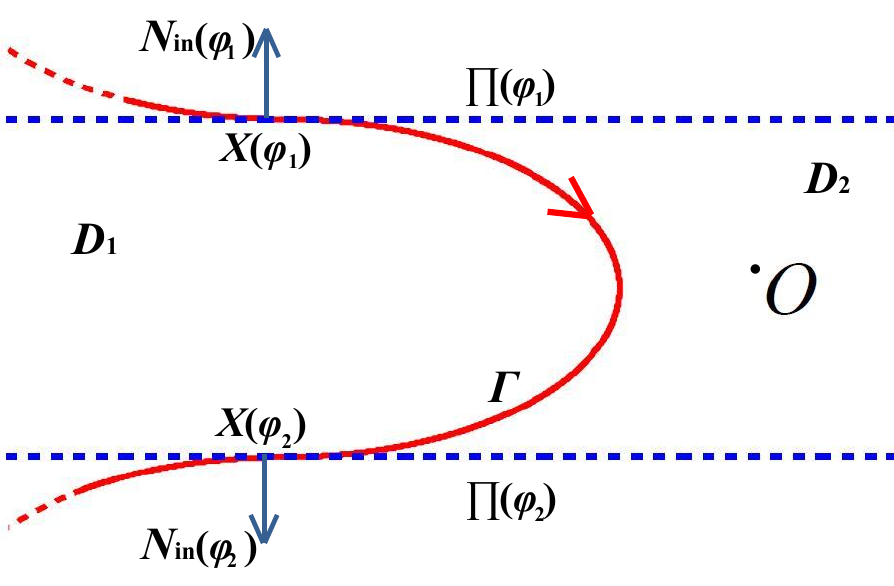}
\caption{An Arc of a ``Star-shaped Curve".}\label{fig:3}
\end{figure}
\begin{lemma}\label{lem:2.4.202401}
Let $X_0$ be an embedded and star-shaped curve. This curve has Chow-Liou-Tsai turning angle greater than $-\pi$.
\end{lemma}
\begin{proof}
Suppose the Chow-Liou-Tsai turning angle of the curve $X_0$ is less than or equal to $-\pi$. There
exists a shortest arc $\Gamma$ of the curve such that the curvature $\kappa \leq 0$ holds everywhere on this arc and
\begin{eqnarray*}
\int_{\Gamma} \kappa ds = -\pi.
\end{eqnarray*}
So the tangent lines at the two end points of $\Gamma$ are parallel.

Let $D$ be the banded region between these two lines.
Denote by $X_0(\varphi_1)$ and $X_0(\varphi_2)$ the two end points of $\Gamma$. The arc $\Gamma$ divides the domain $D$ into two parts. One is convex, denoted by $D_1$, and the other is nonconvex,
denoted by $D_2$. Since the curvature is nonpositive alone the arc $\Gamma$, the star center $O$ lies in the nonconvex domain $D_2$.
However, in this situation, the intersection of half planes $\Pi(\varphi_1)$ and $\Pi(\varphi_2)$ is an empty set. By Lemma \ref{lem:2.2.202401},
this is contradict to the star-shaped property of $X_0$.
\end{proof}

\begin{remark}\label{rem:2.5.202401}
We notice that the conclusion proved in the above lemma was mentioned without proof in the paper by Chow, Liou and Tsai \cite{Chow-Liou-Tsai-1996}. So
Lemma \ref{lem:2.4.202401} essentially belongs to them. The above proofs of Lemma \ref{lem:2.2.202401} and Lemma \ref{lem:2.4.202401} are expository work of their statement.
\end{remark}

Now we turn to prove Theorem \ref{thm:1.1.202401}. Lemma \ref{lem:2.4.202401} says that an embedded and star-shaped
curve has Chow-Liou-Tsai turning angle greater than $-\pi$. So it follows from
Lemma \ref{lem:2.1.202401} that GAPF (\ref{eq:1.1.202401}), with this initial curve, exists globally. Theorem 4.1 in the
paper \cite{Gao-Pan-2023} shows that the curvature of the evolving curve tends to a positive constant
$\sqrt{\frac{\pi}{A_0}}$, where $A_0$ is the area bounded by the initial curve $X_0$. Therefore, there exists a time
$t_0>0$ such that the evolving curve $X(\cdot, t)$ is convex for all $t> t_0$. In the section 4 of the
paper \cite{Gao-Zhang-2019}, we follow Gage's idea \cite{Gage-1986} to prove that the evolving curve converges to a circle as time tends to
infinity. Other details about the convergence of the flow (\ref{eq:1.1.202401}) with a convex initial
curve can be found in the papers \cite{Chao-Ling-Wang-2013} and \cite{Gage-1986}.

\section{A flying wing curve}
In this section, a smooth, embedded and star-shaped curve $X_0$ is constructed to show that Gage's area-preserving flow may not preserve the star shape of the evolving curve.

Let $X(\cdot, t): S^1 \times [0, \omega) \rightarrow \mathbb{R}^2$ be a family of smooth and embedded curve evolving by GAPF as follows
\begin{equation}\label{eq:3.1.202401}
\left\{\begin{array}{ll}
\frac{\partial X}{\partial t}(\varphi, t)=\alpha(\varphi, t) T(\varphi, t)+(\kappa(\varphi, t)-\frac{2\pi}{L(t)})N(\varphi, t) \ \ \ \text{in} \ \ S^1\times (0, \omega),\\
X(\varphi, 0)= X_0(\varphi) \ \  \ \ \ \  \text{on} \ \ S^1,
\end{array} \right.
\end{equation}
where the initial curve $X_0$ is star-shaped; $\alpha(\varphi, t)$ is properly chosen later and $N(\varphi, t)$ stands for the inner unit normal.
According to Proposition 1.1 of the monograph \cite{Chou-Zhu},
the solution to (\ref{eq:3.1.202401}) differs from that of the flow (\ref{eq:1.1.202401}) by a reparametrization. If the curve $X(\cdot, t)$ is star-shaped
with respect to $O$ for $t\in [0, t_0)$, then we may choose
$\alpha=-\frac{\beta}{r}\frac{\partial r}{\partial s}g=-\frac{\beta}{r}\frac{\partial r}{\partial \theta}$ such that the polar angle $\theta$ is independent
of the time $t$ (see Section 2 of the paper \cite{Gao-Pan-2023}). In this case, the radial function $r(\theta, t)$ evolves according to
\begin{eqnarray}\label{eq:3.2.202401}
\frac{\partial r}{\partial t}=\frac{1}{g^2} \frac{\partial^2 r}{\partial \theta^2}-\frac{2}{rg^2}\left(\frac{\partial r}{\partial \theta}\right)^2
-\frac{r}{g^2}+\frac{2\pi g}{rL}, ~~~~(\theta, t)\in S^1 \times [0, \omega).
\end{eqnarray}
In order to explaining the asymptotic behavior of $X(\cdot, t)$, we set $u(\theta, t) = r^2(\theta, t)$. Then we have
$$\frac{\partial u}{\partial \theta} = 2r \frac{\partial r}{\partial \theta},
~~\frac{\partial^2 u}{\partial \theta^2} = 2r \frac{\partial^2 r}{\partial \theta^2} + 2\left(\frac{\partial r}{\partial \theta}\right)^2.$$
So the metric of the evolving curve is
$$g(\theta, t) = \sqrt{r^2 + \left(\frac{\partial r}{\partial \theta}\right)^2} = \sqrt{\frac{4u^2 + \left(\frac{\partial u}{\partial \theta}\right)^2}{4u}}.$$
It follows from the evolution equation of $r(\theta, t)$ that $u(\theta, t)$ satisfies
\begin{eqnarray}\label{eq:3.3.202401}
\frac{\partial u}{\partial t} &=& \frac{4u}{4u^2 + \left(\frac{\partial u}{\partial \theta}\right)^2} \frac{\partial^2 u}{\partial \theta^2}
-\frac{6}{4u^2 + \left(\frac{\partial u}{\partial \theta}\right)^2}\left(\frac{\partial u}{\partial \theta}\right)^2
-\frac{8u^2}{4u^2 + \left(\frac{\partial u}{\partial \theta}\right)^2}
\nonumber\\
&&+\frac{2\pi}{L(t)} \frac{\sqrt{4u^2 + \left(\frac{\partial u}{\partial \theta}\right)^2}}{\sqrt{u}}.
\end{eqnarray}

\begin{figure}[tbh]
\centering
\includegraphics[scale=0.5]{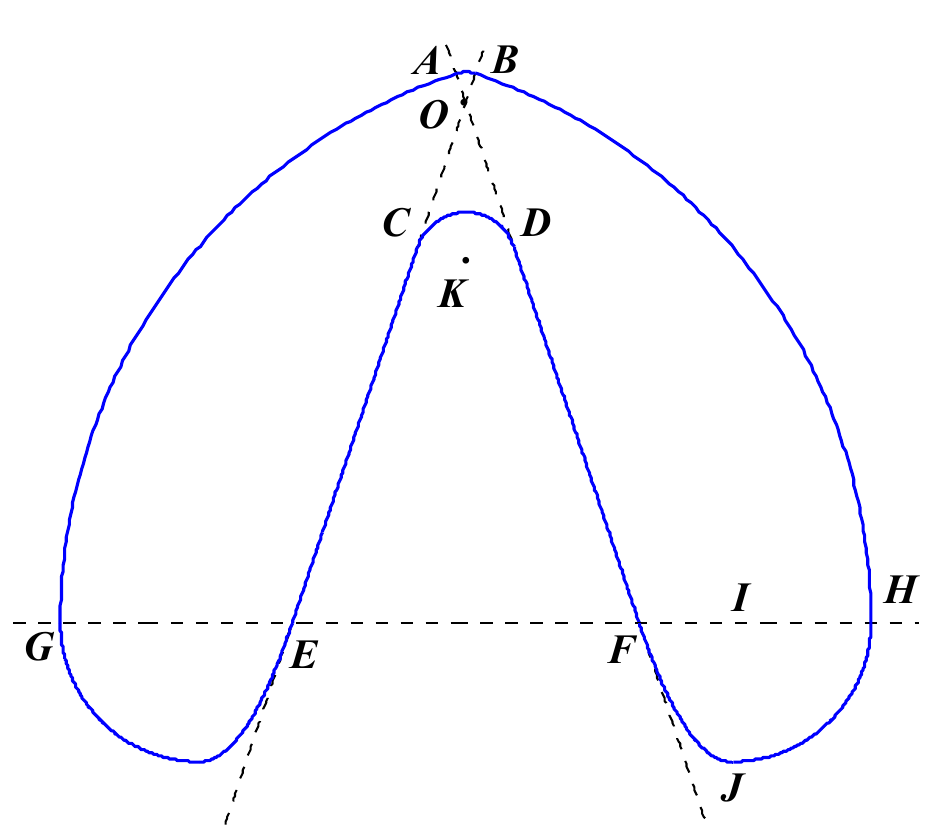}
\caption{A Star-shaped Flying Wing Curve.}\label{fig:4}
\end{figure}
Let $a, b$ be two positive numbers satisfying $2a <b$ and $\rho: = \sqrt{a^2+b^2} <1$. We connect two points $A(-a, b)$ and $B(a, b)$ in the plane by a circular arc
which has radius $\rho$ and center $O$. Denote by $M$ the middle point of this circular arc.
Let $l_1, l_2$ be two positive numbers such that $l_1< l_2$. Let $C(-l_1a, -l_1b)$ and $E(-l_2a, -l_2b)$ be two points on the line
$BO$ and $D(l_1a, -l_1b)$ and $F(l_2a, -l_2b)$ be two points on the line $AO$. Then the length of the segment $EB$ is $|EB|=(l_2+1)\rho$.

Let $H((l_2+1)\rho-l_2a, -l_2b)$ be a point on the line $EF$ so that $|EH|=|EB|$. We connect $B$ and $H$ by a circular arc with radius $|EH|=(l_2+1)\rho$ and center $E$.
Let $I$ be a point on the segment $FH$ with $x$-coordinate $x_I = (1-\lambda)l_2 a + \lambda [(l_2+1)\rho -l_2a]$, where $\lambda \in (0, 1)$ is to be determined later. So the length of $IH$ is
$|IH|=(1 - \lambda)[(l_2+1)\rho -2l_2a]$. Let $J$ be the point on the vertical line passing through the point $I$ and $|IJ|=|IH|$. Then the coordinate of $J$ is given by
\begin{eqnarray}\label{eq:3.4.202401}
(x_J, y_J) = \big((1-\lambda)l_2a+\lambda [(l_2+1)\rho-l_2a], -l_2b-(1-\lambda)[(l_2+1)\rho-2l_2a] \big).
\end{eqnarray}
We connect two points $H$ and $J$ in the plane by a circular arc which has radius $|IH|$ and center $I$.

Then we connect $F$ and $J$ by a parabola arc which has vertex $J$, symmetric line $IJ$ and is tangent to the line $DF$ at the point $F$.
We set the equation of this parabola arc as
\begin{eqnarray}\label{eq:3.5.202401}
y = k \big[x- (1-\lambda)l_2a- \lambda[(l_2+1)\rho-l_2a] \big]^2 - l_2b - (1-\lambda)[(l_2+1)\rho-2l_2a],
\end{eqnarray}
where the constant $k$ is positive and $x\in [x_F, x_I]$. This parabola arc passes through the point $F$, so
$$k \lambda^2 \big[(l_2+1)\rho-2l_2a \big]^2 = (1-\lambda)[(l_2+1)\rho-2l_2a],$$
i.e.,
\begin{eqnarray}\label{eq:3.6.202401}
k = \frac{1- \lambda}{\lambda^2} \frac{1}{(l_2+1)\rho-2l_2a} >0.
\end{eqnarray}
Since the parabola arc is also tangent to the segment $DF$ at the point $F$, we have
$$2k\big[x_F- (1-\lambda)l_2a- \lambda [(l_2+1)\rho-l_2a] \big] = -\frac{b}{a},$$
which, together with \eqref{eq:3.6.202401} and the fact $x_F = l_2a$, implies that
\begin{eqnarray}\label{eq:3.7.202401}
\lambda = \frac{2a}{2a +b}.
\end{eqnarray}
The numbers $a$ and $b$ are fixed. If we set $\lambda$ the value given by \eqref{eq:3.7.202401} and $k$ given in \eqref{eq:3.6.202401}, then we obtain
a parabola arc $JF$ which is tangent to the segment $DF$ at the point $F$.

We connect $F$ and $D$ by a segment.
Let $BP$ be the segment perpendicular to the segment $OM$ at the point $P$, where $M$ is the middle point of the circular arc $AB$.
Then $|BP|= a$ and $|OP|= b$. There is a unique line passing through
$D$ and perpendicular to $OF$ at the point $D$. There is also a unique line passing through $C$ and perpendicular to $OE$ at the point $C$. The two
line intersect at a point, denoted by $K$. So there exists a circular arc from $C$ to $D$ with radius $|KD|$ and center $K$. Noticing that the right triangle
Rt$\triangle OBP$ is similar to the right triangle Rt$\triangle OKD$, we have $\frac{|OK|}{|OB|} = \frac{|OD|}{|OP|}$, i.e.,
$$|OK| = \frac{l_1 \rho^2}{b}.$$
Let $N$ be the middle point of the circular arc $CD$. Now we have a $C^1$ curve $MBHJFDN$. Together with its symmetric image curve with respect to $y$-axis,
we have constructed a $C^1$, closed and embedded curve, denoted by $Y_0$. It follows from Lemma \ref{lem:2.2.202401} that this curve is star-shaped.
The set of its star centers is the open domain bounded by the segments $OA$ and $OB$ and the circular arc $AB$.
The domain bounded by $Y_0$ is like a flying wing. An example is shown in Figure \ref{fig:4} with parameters $a= 0.1, b=0.3, l_1=5, l_2=18$. The name "flying wing"
curve is motivated by Lai's recent work \cite{Lai-2024}.

If we evolve the $C^1$ curve $Y_0$ as an initial curve by the curve shortening flow, then it follows from the results in the paper \cite{Grayson-1987} that
there is a family of smooth and embedded curve $Y(\cdot, t)$ for $t\in (0, \omega)$. Let $A_t$ and $B_t$ be the points on the curve $Y(\cdot, t)$
with initial values $A$ and $B$, respectively. By the continuity, the curve $Y(\cdot, t)$ converges to
$Y_0$ point wise as $t \rightarrow 0$ and $Y(\cdot, t)$ is star-shaped for small $t>0$. So there exists a very small $\varepsilon>0$ such that
the curve $X_0(\cdot):= Y(\cdot, \varepsilon)$ is smooth and star-shaped with respect to a point, which is chosen as the origin of the plane from now on.

\begin{example}\label{exa:3.1.202401}
Consider Gage's area-preserving flow with initial $X_0$ defined by above $Y(\cdot, \varepsilon)$ for some very small $\varepsilon>0$. By the continuity,
we have on the arc $A_\varepsilon B_\varepsilon$,
$$u(\theta, 0) \approx \rho^2, ~~\frac{\partial u}{\partial \theta}(\theta, 0) \approx 0, ~~\frac{\partial^2 u}{\partial \theta^2}(\theta, 0) \approx 0.$$
So the right hand side of the equation \eqref{eq:3.3.202401} can be approximated by $-2 + \frac{4\pi}{L(0)}\rho$ on the arc $A_\varepsilon B_\varepsilon$.
Now we may pick proper values of $a, b, l_1$ and $l_2$ such that the initial length $L(0)$ is large and $\rho$ is very small and $\frac{2\pi}{L(0)}\rho$
is sufficiently closed to 0. Hence, under Gage's area-preserving flow \eqref{eq:3.1.202401}, the function
$u_{\min}(t): = \min\{u(\theta, t)| \theta \in S^1\}$, with quite small initial value $u_{\min}(0)\approx \rho^2$, may drop to 0 in a short time.
In the same time, the arc near $A_\varepsilon B_\varepsilon$ may sweep through the kernel of the evolving curve, which makes the evolving curve no longer star-shaped.
\end{example}

In the above analysis, a star-shaped flying wing curve $X_0$ is constructed as an initial curve for Gage's area-preserving flow which may lose the star shape of the evolving curve $X(\cdot, t)$.
Next, we do some more analysis on the kernel of the evolving curve along GAPF \eqref{eq:3.1.202401}. And we will show that
the evolving curve indeed loses its star shape in a short time.

Denote by $Z(\theta, t, \cdot)$ the tangent line of the curve $X(\cdot, t)$ at the point $(\theta, t)$. The expression of $Z$ is given by
\begin{eqnarray*}
Z(\theta, t, \lambda) = X(\theta, t) + \lambda T(\theta, t), ~~\theta\in S^1,~~\lambda\in \mathbb{R}.
\end{eqnarray*}
Along GAPF \eqref{eq:3.1.202401}, the tangent line $Z(\theta, t, \cdot)$ moves according to
\begin{eqnarray*}
\frac{\partial Z}{\partial t}(\theta, t, \lambda) &=& \frac{\partial X}{\partial t}(\theta, t) + \lambda \frac{\partial T}{\partial t}(\theta, t)
\\
&=& - \frac{1}{r} \frac{\partial r}{\partial \theta}\beta T +\beta N
+ \lambda \left(- \frac{1}{r} \frac{\partial r}{\partial \theta}\beta \kappa + \frac{\partial \kappa}{\partial \theta} \frac{1}{g}\right)  N,
\end{eqnarray*}
where $\beta (\theta, t)= \kappa(\theta, t) - \frac{2\pi}{L(t)}$.
The boundary of the kernel $\mathfrak{C}(t)$ consists of arcs of the curve $X(\cdot, t)$ and segments of the tangent lines.
In the later cases, the curvature at the tangent point $(\theta, t)$ is 0. Alone the segments $A_{\varepsilon}O $ and $OB_{\varepsilon}$
in the Example \ref{exa:3.1.202401}, we have $\kappa(\theta, t)=0$ at the tangent points.
And alone the arc $\overset{\frown}{B_{\varepsilon}A_{\varepsilon}}$, we have $\lambda=0$. This observation belongs to Halpern \cite{Halpern-1969}.

Let $\widetilde{A}(t)$ be the area of the kernel $\mathfrak{C}(t)$. The first variation formula implies that the area $\widetilde{A}(t)$ varies according to
\begin{eqnarray*}
\frac{d \widetilde{A}}{d t}(t) &=& -\int_{\overset{\frown}{B_{\varepsilon}A_{\varepsilon}}} \left(\kappa - \frac{2\pi}{L(t)} \right)
- \int_{A_{\varepsilon}O \cup OB_{\varepsilon}} \left \langle \frac{\partial Z}{\partial t}, N\right\rangle
\\
&=& -\int_{\overset{\frown}{B_{\varepsilon}A_{\varepsilon}}} \left(\kappa - \frac{2\pi}{L(t)} \right)ds
- \int_{A_{\varepsilon}O \cup OB_{\varepsilon}} \left(- \frac{2\pi}{L(t)} + \lambda\frac{\partial \kappa}{\partial \theta} \frac{1}{g} \right) du
\\
&\leq& -\int_{\overset{\frown}{B_{\varepsilon}A_{\varepsilon}}} \kappa ds + \frac{2\pi}{L(t)} |\overset{\frown}{B_{\varepsilon}A_{\varepsilon}}|
+\frac{4\pi}{L(t)} |A_{\varepsilon}O| + \left|\frac{\partial \kappa}{\partial \theta}\right| \frac{2}{g} |A_{\varepsilon}O|.
\end{eqnarray*}
Since $\int_{\overset{\frown}{B_{\varepsilon}A_{\varepsilon}}} \kappa ds \rightarrow 2 \arctan \frac{a}{b}$ as $\varepsilon \rightarrow 0^+$, we may choose small $\varepsilon>0$
so that, for small $t>0$ next estimate holds:
\begin{eqnarray*}
\frac{d \widetilde{A}}{d t}(t) &\leq& -\frac{3}{2} \arctan\frac{a}{b}+\frac{\pi}{L(t)}\rho^2 \arctan\frac{a}{b}
+ \frac{4\pi}{L(t)} \rho + \left|\frac{\partial \kappa}{\partial \theta}\right| \frac{2}{g} \rho.
\end{eqnarray*}
The evolving curve is smooth, so one may choose very small $a$ and $b$ and very large $l_1$ and $l_2$ so that the last three terms in the above inequality
are small enough comparing to a term $\frac{1}{2}\arctan\frac{a}{b}$.

For example, one may choose $a \leq 10^{-10}$ and $b=4a$. So $\rho = \sqrt{a^2 +b^2} < 10^{-9}$. Let $l_1$ be large, e.g., $l_1 = 10^{10}$
and $l_2$ be larger, e.g. $l_2 \geq 10^{100}$. Then $\widetilde{A}(0) \leq 10^{-9}$ is small and for small $t>0$, the derivative of $\widetilde{A}(t)$ satisfies
\begin{eqnarray*}
\frac{d \widetilde{A}}{d t}(t) < - \arctan\frac{1}{4} = - 0.24497\cdots.
\end{eqnarray*}
In this setting, the area $\widetilde{A}(t)$ of the kernel $\mathfrak{C}(t)$ drops to 0 in a short time interval.
As we know, the kernel of a star-shaped curve is a convex set. So the closure $\overline{\mathfrak{C}(t)}$ will degenerate to a segment or a point.
Therefore, along GAPF, the curve $X(\cdot, t)$ will lose its star shape, i.e., there is $t_0>0$ so that $X(\cdot, t)$ is not star-shaped with respect to any points
for $t>t_0$.

\section{Star-shaped curves under the CSF}

In the paper \cite{Gao-Pan-2023}, the authors study centrosymmetric and star-shaped curves evolving according to Gage's area-preserving flow,
by comparing the behavior of the CSF. Next, we follow the similar idea to study star-shaped curves evolving according to
the CSF.

Let $X_0$ be an embedded and star-shaped curve in the plane with a star center $O$.
We evolve $X_0$ as the initial curve according to the CSF. Then we obtain a family of smooth and
embedded curve $\widetilde{X}(\cdot, t)$ and the evolving curve is star-shaped for small $t>0$. We revise the curve shortening flow as follows:
\begin{equation}\label{eq:4.1.202401}
\left\{\begin{array}{ll}
\frac{\partial \widetilde{X}}{\partial t}(\varphi, t)=-\frac{\widetilde{\kappa}}{\widetilde{r}}\frac{\partial \widetilde{r}}{\partial \theta} \widetilde{T}(\varphi, t)
  +\widetilde{\kappa}(\varphi, t)\widetilde{N}(\varphi, t)   \ \ \ \text{in} \ \ S^1\times (0, \omega),\\
\widetilde{X}(\varphi, 0)= X_0(\varphi) \ \  \ \ \ \  \text{on} \ \ S^1.
\end{array} \right.
\end{equation}
Under the flow \eqref{eq:4.1.202401}, the polar angle $\theta$ is independent of the time. And the radial function $\widetilde{r}(\theta, t)$ evolves according to
\begin{eqnarray}\label{eq:4.2.202401}
\frac{\partial \widetilde{r}}{\partial t}=\frac{1}{\widetilde{g}^2} \frac{\partial^2 \widetilde{r}}{\partial \theta^2}
-\frac{2}{\widetilde{r}\widetilde{g}^2}\left(\frac{\partial \widetilde{r}}{\partial \theta}\right)^2 -\frac{\widetilde{r}}{\widetilde{g}^2}, ~~\text{on}~ S^1\times [0, \omega).
\end{eqnarray}

Let $X_0$ evolve as an initial curve of GAPF (\ref{eq:1.1.202401}) and the CSF \eqref{eq:4.1.202401}, respectively.
We obtain two family of smooth and embedded curves, denoted by $X(\cdot, t)$ and $\widetilde{X}(\cdot, t)$, where $t\in [0, \omega)$, $\omega= \frac{A_0}{2\pi}$ and $A_0$
is the area bounded by $X_0$.
\begin{lemma}\label{lem:4.1.202401}
For each $t\in (0, \omega)$, the curve $\widetilde{X}(\cdot, t)$ is strictly enclosed by the curve $X(\cdot, t)$,
i.e., $\widetilde{X}(\cdot, t)$ lies in the open domain bounded by the curve $X(\cdot, t)$.
\end{lemma}
\begin{proof}
By the continuity, there is a small positive $t_1$ such that both $X(\cdot, t)$ and $\widetilde{X}(\cdot, t)$ are star-shaped with respect to $O$ in the time interval $[0, t_1)$.
Under the two flows, the radial functions $r(\theta, t)$ and $\widetilde{r}(\theta, t)$ satisfy the equations \eqref{eq:3.2.202401} and \eqref{eq:4.2.202401}, separately.
Set $f(\theta, t) = r(\theta, t) - \widetilde{r}(\theta, t)$. Then this function satisfies (see also the equation (3.16) in the paper \cite{Gao-Pan-2023})
\begin{eqnarray}
\frac{\partial f}{\partial t}
&=& \frac{1}{g^2}\frac{\partial^2 f}{\partial \theta^2}
-\frac{1}{\widetilde{g}^2 g^2}\frac{\partial^2 \widetilde{r}}{\partial \theta^2}
           \left(\frac{\partial \widetilde{r}}{\partial \theta} + \frac{\partial r}{\partial \theta}\right)
           \frac{\partial f}{\partial \theta}
-\frac{2}{\widetilde{r} \widetilde{g}^2}\left(\frac{\partial \widetilde{r}}{\partial \theta}
+\frac{\partial r}{\partial \theta}\right)\frac{\partial f}{\partial \theta}
\nonumber\\
&& +\frac{2}{\widetilde{r} \widetilde{g}^2 g^2}\left(\frac{\partial \widetilde{r}}{\partial \theta}\right)^2 \left(\frac{\partial \widetilde{r}}{\partial \theta}
+\frac{\partial r}{\partial \theta}\right)\frac{\partial f}{\partial \theta}
+\frac{r}{\widetilde{g}^2 g^2}
    \left(\frac{\partial \widetilde{r}}{\partial \theta} + \frac{\partial r}{\partial \theta}\right)\frac{\partial f}{\partial \theta}
\nonumber\\
&& -\frac{r + \widetilde{r}}{\widetilde{g}^2 g^2} \frac{\partial^2 \widetilde{r}}{\partial \theta^2} f
+\frac{2(r + \widetilde{r})}{\widetilde{r} \widetilde{g}^2 g^2} \left(\frac{\partial \widetilde{r}}{\partial \theta}\right)^2 f
+\frac{2}{\widetilde{r} r g^2} \left(\frac{\partial \widetilde{r}}{\partial \theta}\right)^2 f
\nonumber\\
&&
-\frac{f}{g^2}
+\frac{r(r + \widetilde{r})}{\widetilde{g}^2 g^2}f
+\frac{2\pi g}{rL}.
\label{eq:4.3.202401}
\end{eqnarray}
This is a parabolic equation with smooth and bounded coefficients on the domain $S^1 \times [0, t_1)$, where $0< t_1 < \omega$.
For $0<t\leq t_1$, suppose $f(\theta, t)$ attains its minimum w.r.t. $\theta$ at a point $\theta_*$. Then $\frac{\partial^2 f}{\partial \theta^2}(\theta_*, t) \leq 0$,
and $\frac{\partial f}{\partial \theta}(\theta_*, t) =0$. Since both $X(\cdot, t)$ and $\widetilde{X}(\cdot, t)$ are smooth and star-shaped w.r.t. $O$ at the
time interval $[0, t_1]$, there exists a positive number $C=C(t_1)$ so that at the point $(\theta_*, t)$ we have
\begin{eqnarray*}
\frac{\partial f}{\partial t}(\theta_*, t) \geq -C f(\theta_*, t)+\sqrt{\frac{\pi}{A_0}}.
\end{eqnarray*}
Denote $f_{\min}(t)$ by the minimum value of $f(\theta, t)$ with respect to $\theta$. Then $f_{\min}(0)=0$ and it follows from the maximum principle,
\begin{eqnarray*}
f_{\min}(t) \geq \frac{1}{C}\sqrt{\frac{\pi}{A_0}}(1-e^{-Ct}).
\end{eqnarray*}
So the function $f$ is positive everywhere for $t>0$. And we have on the time interval $t\in (0, t_1)$ next estimate
\begin{eqnarray} \label{eq:4.4.202401}
r (\theta, t) > \widetilde{r}(\theta, t)>0, ~~\theta\in S^1.
\end{eqnarray}
This means that the curve $\widetilde{X}(\cdot, t)$ is contained in the open domain enclosed by the curve $X(\cdot, t)$ for every small $t>0$.

Next we show that $X(\cdot, t)$ and $\widetilde{X}(\cdot, t)$ never intersect for $t>0$ until the latter shrinks to a point.
Suppose there exists a smallest time $t_2>0$ such that $X(\cdot, t_2)$ and $\widetilde{X}(\cdot, t_2)$ intersect at a tangent point.
Near this point, we locally express two curves as two graphs $y = y(x, t)$ and $\widetilde{y} = \widetilde{y}(x, t)$ on $D:=[-a, a]\times [t_2- \delta, t_2]$.
The direction of $x$-axis is $T$ and the direction of $y$-axis is the inner normal $N$.
Set $w(x, t) = y(x, t) - \widetilde{y}(x, t)$. Without loss of generality, we assume $w(x, t) \leq 0$ on the domain $D$. There is a point $x_*\in (a, b)$ such that $w(x_*, t_2) =0$.
We now deduce the evolution equation of the function $w(x, t)$.

Suppose we have a family of smooth curves $Y(\cdot, t)$ satisfying next equation
\begin{eqnarray}\label{eq:4.5.202401}
\frac{\partial Y}{\partial t}(x, t) =  \alpha(x, t) T(x, t) + \beta(x, t)  N(x, t)
     \ \ \textup{in}~ (x_0-\delta, x_0+\delta) \times (0, t_0) : = D,
\end{eqnarray}
where $Y(x, t) = (x, y(x, t))$; the Frenet frame is as follows
\begin{eqnarray*}
T(x, t) = \frac{1}{\sqrt{1+\left(\frac{\partial y}{\partial x}\right)^2}}\left(1, \frac{\partial y}{\partial x}\right), ~~
N(x, t) = \frac{1}{\sqrt{1+\left(\frac{\partial y}{\partial x}\right)^2}}\left(-\frac{\partial y}{\partial x}, 1\right);
\end{eqnarray*}
and the curvature is given by
\begin{eqnarray*}
\kappa(x, t) = \frac{1}{\left(1+\left(\frac{\partial y}{\partial x}\right)^2\right)^{3/2}} \frac{\partial^2 y}{\partial x^2}.
\end{eqnarray*}
So the equation \eqref{eq:4.5.202401} can be written as
\begin{eqnarray*}
\left(\frac{\partial x}{\partial t}, \frac{\partial y}{\partial t}\right) &=& \frac{\alpha}{g} \left(1, \frac{\partial y}{\partial x}\right)
+ \frac{\beta}{g} \left(-\frac{\partial y}{\partial x}, 1\right)
= \frac{1}{g} \left(\alpha - \beta \frac{\partial y}{\partial x}, \alpha\frac{\partial y}{\partial x} + \beta \right).
\end{eqnarray*}
This gives us
\begin{eqnarray*}
\frac{\partial x}{\partial t}=\frac{1}{g} \left(\alpha - \beta \frac{\partial y}{\partial x}\right).
\end{eqnarray*}
If we choose $\alpha = \beta \frac{\partial y}{\partial x}$ then $\frac{\partial x}{\partial t} \equiv 0$ and the evolution equation of $y(x, t)$ is
\begin{eqnarray}\label{eq:4.6.202401}
\frac{\partial y}{\partial t} = \frac{1}{g} \left(\alpha\frac{\partial y}{\partial x} + \beta \right) = \beta g.
\end{eqnarray}

If $\beta = \kappa - \frac{2\pi}{L}$ then, under  GAPF, $y(x, t)$ satisfies
\begin{eqnarray}\label{eq:4.7.202401}
\frac{\partial y}{\partial t} = \frac{1}{1+\left(\frac{\partial y}{\partial x}\right)^2} \frac{\partial^2 y}{\partial x^2} - \frac{2\pi}{L(t)}g.
\end{eqnarray}
If $\beta = \kappa$ then, under the CSF, $\widetilde{y}(x, t)$ satisfies
\begin{eqnarray}\label{eq:4.8.202401}
\frac{\partial \widetilde{y}}{\partial t} = \frac{1}{1+\left(\frac{\partial \widetilde{y}}{\partial x}\right)^2} \frac{\partial^2 \widetilde{y}}{\partial x^2}.
\end{eqnarray}
Therefore, it follows from \eqref{eq:4.7.202401} and \eqref{eq:4.8.202401} that the function $w=y- \widetilde{y}$ satisfies next evolution equation on the domain $D$:
\begin{eqnarray}\label{eq:4.9.202401}
\frac{\partial w}{\partial t} = \frac{1}{1+\left(\frac{\partial \widetilde{y}}{\partial x}\right)^2} \frac{\partial^2 w}{\partial x^2}
- \frac{\frac{\partial^2 y}{\partial x^2}\left(\frac{\partial y}{\partial x}+\frac{\partial \widetilde{y}}{\partial x}\right)}
{\left(1+\left(\frac{\partial y}{\partial x}\right)^2\right)\left(1+\left(\frac{\partial \widetilde{y}}{\partial x}\right)^2\right)} \frac{\partial w}{\partial x}
- \frac{2\pi}{L(t)}g.
\end{eqnarray}
Since $w(x, t)\leq -\varepsilon_0 <0$ holds on the parabolic boundary
$$\big([-a, a]\times \{t_2- \delta\}\big) \cup \big(\{-a\}\times [t_2- \delta, t_2]\big) \cup \big(\{a\} \times [t_2- \delta, t_2]\big)$$
for some constant $\varepsilon_0 >0$, the maximum principle implies that $w(x, t) \leq -\varepsilon_0 <0$ holds on the whole domain $D$.
This is in contradiction with the assumption $w(x_*, t_2) =0$ for some $x_*\in (-a, a)$.
\end{proof}

\begin{lemma}\label{lem:4.2.202401}
Let $X_0$ be an embedded and star-shaped curve in the plane. Let this curve evolve as an initial curve of Gage's area-preserving flow (\ref{eq:1.1.202401})
and also the curve shortening flow \eqref{eq:4.1.202401} separately . We have two family of smooth curves denoted by $X(\cdot, t)$ and $\widetilde{X}(\cdot, t)$, respectively.
If $X(\cdot, t)$ and $\widetilde{X}(\cdot, t)$ are star-shaped curves on the same time interval $[0, t_*)$, then
the kernel of $\widetilde{X}(\cdot, t)$ is a subset of the kernel of the curve $X(\cdot, t)$, i.e.,
\begin{eqnarray}\label{eq:4.10.202401}
\mathfrak{C}(\widetilde{X}(\cdot, t))  \subseteq \mathfrak{C}(X(\cdot, t)),~~ t\in [0, t_*), ~~t_*>0.
\end{eqnarray}
\end{lemma}
\begin{proof}
For $t=0$, the set $\mathfrak{C}(\widetilde{X}(\cdot, t))$ coincides with  $\mathfrak{C}(X(\cdot, t))$.
We first show that \eqref{eq:4.10.202401} holds for each small $t>0$.
The boundary of the kernel of the a star-shaped curve consists of segments and arcs of the form
\begin{eqnarray}\label{eq:4.11.202401}
Z(\varphi, t) = X(\varphi, t) + \lambda T(\varphi, t), ~~\lambda\in [a_\varphi, b_\varphi].
\end{eqnarray}
If the boundary point $Z(\varphi, t)$ lies on the curve $X(\cdot, t)$, then we have $\lambda =0$. Under the flow
$\frac{\partial X}{\partial t}(\varphi, t) =  \beta(\varphi, t)  N(\varphi, t)$,
we have
\begin{eqnarray}\label{eq:4.12.202401}
\frac{\partial Z}{\partial t}(\varphi, t) = \left(\beta(\varphi, t) + \lambda \frac{\partial \beta}{\partial s}(\varphi, t) \right)  N(\varphi, t).
\end{eqnarray}
Therefore the boundary of $\mathfrak{C}(X(\cdot, t))$ under Gage's area-preserving flow moves according to
\begin{eqnarray}\label{eq:4.13.202401}
\frac{\partial Z}{\partial t}(\varphi, t) = \left(\kappa(\varphi, t) - \frac{2\pi}{L(t)}+ \lambda \frac{\partial \kappa}{\partial s}(\varphi, t) \right)  N(\varphi, t),
\end{eqnarray}
and boundary of $\mathfrak{C}(\widetilde{X}(\cdot, t))$ under the curve shortening flow satisfies
\begin{eqnarray}\label{eq:4.14.202401}
\frac{\partial \widetilde{Z}}{\partial t}(\varphi, t) = \left(\widetilde{\kappa}(\varphi, t)
  + \lambda \frac{\partial \widetilde{\kappa}}{\partial \widetilde{s}}(\varphi, t) \right)  \widetilde{N}(\varphi, t).
\end{eqnarray}
Since, at time $t=0$, we have
\begin{eqnarray}\label{eq:4.15.202401}
\frac{\partial Z}{\partial t}(\varphi, 0) - \frac{\partial \widetilde{Z}}{\partial t}(\varphi, 0) = - \frac{2\pi}{L(0)}  N(\varphi, 0),
\end{eqnarray}
it follows from the first variation formula that there exists $t_3>0$ such that \eqref{eq:4.10.202401} holds on time interval $[0, t_3)$.

Suppose there exists a smallest time $t_4>0$ such that $\mathfrak{C}(\widetilde{X}(\cdot, t_4)) \nsubseteq \mathfrak{C}(X(\cdot, t_4))$, i.e., there is point
$P \in \mathfrak{C}(\widetilde{X}(\cdot, t_4))$ but $P\notin \mathfrak{C}(X(\cdot, t_4))$. $\mathfrak{C}(\widetilde{X}(\cdot, t_4))$ is an open set, so
there is $\eta>0$ such that $P \in \mathfrak{C}(\widetilde{X}(\cdot, t))$ holds for $t\in [t_4-\eta, t_4]$. By the smallest property assumption of $t_4>0$,
there exists $\varepsilon>0$ such that $P \in \mathfrak{C}(X(\cdot, t))$ holds for $t\in [t_4-\varepsilon, t_4)$.
Take $\xi=\min\{\eta, \varepsilon\}$. We have $P \in \mathfrak{C}(X(\cdot, t))\cap \mathfrak{C}(\widetilde{X}(\cdot, t))$ holds for $t\in [t_4-\xi, t_4)$.

Choose $P$ as the original point of the plane. By Lemma \ref{lem:4.1.202401}, $X(\cdot, t)$ contains $\widetilde{X}(\cdot, t)$ as $t\in [t_4-\xi, t_4)$. This implies that
$r(\theta, t_4-\xi) > \widetilde{r}(\theta, t_4-\xi)$. It follows from the proof of Lemma \ref{lem:4.1.202401} we have
\begin{eqnarray}\label{eq:4.16.202401}
r(\theta, t) > \widetilde{r}(\theta, t), ~~(\theta, t)\in S^1 \times [t_4-\xi, t_4).
\end{eqnarray}
By the assumption $P\in \mathfrak{C}(\widetilde{X}(\cdot, t_4))$, $\widetilde{r}(\theta, t_4-\xi)$ has a positive lower bound $\delta>0$. By the
continuity, the equation \eqref{eq:4.16.202401} implies
\begin{eqnarray}\label{eq:4.17.202401}
r(\theta, t_4) > \delta >0, ~~\theta\in S^1.
\end{eqnarray}
This means that the curve $X(\cdot, t)$ does not touch the point $P$ during the time $t\in [t_4-\xi, t_4]$. By the gradient estimate in Lemma 2.14 of the paper
\cite{Gao-Pan-2023}, $P$ is a star center of $X(\cdot, t_4)$. This contradicts the assumption $P\notin \mathfrak{C}(X(\cdot, t_4))$.
\end{proof}

Now we prove Theorem \ref{thm:1.3.202401} as follows. Let a star-shaped curve $X_0$ be an initial curve of GAPF (\ref{eq:1.1.202401}) and
the CSF \eqref{eq:4.1.202401}, respectively. We obtain two family of smooth and embedded curves $X(\cdot, t)$ and $\widetilde{X}(\cdot, t)$,
where $t\in [0, \omega)$ and $\omega = \frac{A_0}{2\pi}$. Whenever $X(\cdot, t)$ and $\widetilde{X}(\cdot, t)$ are star-shaped, Lemma \ref{lem:4.2.202401} says
the kernel of $\widetilde{X}(\cdot, t)$ is a subset of the kernel of the curve $X(\cdot, t)$ for each $t\in (0, \omega)$.

Suppose $X(\cdot, t)$ is a star-shaped curve for $t\in [0, t_*)$ but it is not star-shaped with respect to any points at time $t_*\in (0, \frac{A_0}{2\pi})$. Then the closure of
its convex kernel $\overline{\mathfrak{C}(X(\cdot, t))}$ degenerates to a segment or a point as $t\rightarrow t_*$. If $\widetilde{X}(\cdot, t)$ is not star-shaped with respect to any points
before the time $t_*$, then we have done. Otherwise, by the fact $\mathfrak{C}(\widetilde{X}(\cdot, t))  \subseteq \mathfrak{C}(X(\cdot, t))$,
$\overline{\mathfrak{C}(\widetilde{X}(\cdot, t))}$ also degenerates as a segment or a point as $t \rightarrow t_*$. In this case,
$\widetilde{X}(\cdot, t_*)$ is not a star-shaped curve any more. This completes the proof of Theorem \ref{thm:1.3.202401}.

In conclusion, we show that Gage's area-preserving flow evolves smooth, embedded and star-shaped curves globally on time interval
$[0, +\infty)$. Although in some extreme cases, this flow may not preserve the star shape of the evolving curve.
On the other hand, once GAPF
for smooth, embedded and star-shaped curves exist globally, by the previous studies \cite{Dittberner-2021, Gao-Pan-2023, Gao-Zhang-2019},
the evolving curves converge to circles as $t \rightarrow +\infty$.

\begin{remark}\label{rem:4.3.202401}
The smooth star-shaped curve $X_0$ defined in the Example \ref{exa:3.1.202401} is symmetric with respect to the $y$-axis. If we rotate this
curve with respect to $y$-axis, then we have a smooth star-shaped surface, denoted by $Y_0$, in the 3-dimensional Euclidean space.
It seems that the mean curvature flow with initial surface $Y_0$ may not always preserve star shape of the evolving surface. More analysis
on the evolution of star-shaped surfaces along the mean curvature flow can be found in Mantegazza's note \cite{Mantegazza-2010}.
\end{remark}

~\\
\textbf{Acknowledgments}
During past ten years, Laiyuan Gao persisted in studying GAPF for star-shaped curves, encouraged continuously by his Ph.D. advisor Professor Shengliang Pan.
Gao thanks his mentor Professor Lei Ni for constant concern these years.
As the first author mentioned that to Professor Bonnett Chow,
Gao thanks Professor Dong-Ho Tsai for teaching him a lot via their past collaborations.
Gao also thanks Professor Jianbo Fang for telling him Halpern's result \cite{Halpern-1969}.

{\bf Laiyuan Gao}

School of Mathematics and Statistics, Jiangsu Normal University.

No.101, Shanghai Road, Xuzhou City, Jiangsu Province, China.

Email: lygao@jsnu.edu.cn\\

{\bf Shicheng Zhang}

School of Mathematics and Statistics, Jiangsu Normal University.

No.101, Shanghai Road, Xuzhou City, Jiangsu Province, China.

Email: zhangshicheng@jsnu.edu.cn\\

{\bf Yuntao Zhang}

School of Mathematics and Statistics, Jiangsu Normal University.

No.101, Shanghai Road, Xuzhou City, Jiangsu Province, China.

Email: yuntaozhang@jsnu.edu.cn\\

\end{document}